\documentclass[12pt,twoside,a4paper]{article}
\usepackage[cp866]{inputenc}
\usepackage[russian]{babel}
\usepackage{amsmath,amsfonts,amscd,amssymb,latexsym}

\begin{document}
\newcommand\Der{{\rm Der}}
\newcommand\DD{{\rm D}}
\newcommand\party{\frac{\partial}{\partial y}}
\newcommand\partx{\frac{\partial}{\partial x}}

\vspace*{1cm}

\begin{center}
{\large \bf
New Examples of Simple Jordan Superalgebras\\
over an~Arbitrary Field of Characteristic Zero}
\end{center}

\vspace{7mm}
\begin{center}
{\bf V.~N.~Zhelyabin}
\end{center}

\footnotetext{RFBR 09-01-00157,
SB RAS grant,
Development of Scientific Potential of Higher School
(Grant 2.1.1.419),
Federal Program
``Scientific and Pedagogical Staff of Innovative Russia'' 2009-2013
(State Contract No.~02.740.11.0429). }

\renewcommand{\abstractname}{}
\begin{abstract}
{\bf Abstract:} {\it An~new example of a~unital simple special
Jordan superalgebra over the field of real numbers
 was  constructed in \cite{ShestZhel}.
 It turned out
to be a~subsuperalgebra of the Jordan superalgebra of vector type
 $J(\Gamma,D)$,
but not isomorphic to a~superalgebra of this type. Moreover, its
superalgebra of fractions is isomorphic to a~Jordan superalgebra
of vector type. A~similar example of a~Jordan superalgebra over
a~field of characteristic~0 in which the equation $t^2+1=0$ has no
solutions was constructed in \cite{Zhel09}. In this article we
present an~example of a~Jordan superalgebra with the same
properties over an~arbitrary field of characteristic~0. A~similar
example of a~superalgebra is found in the Cheng--Kac
superalgebra.}
\bigskip

\end{abstract}

\vspace{6pt}

Jordan algebras and superalgebras constitute
an~important class of algebras in ring theory.
Simple Jordan superalgebras are studied in 
\cite{Kac,Kantor,ShAlt,Z2000, MZ2001,KMZ2001, Zel&Rac, CantarKac}.

The unital simple special Jordan superalgebras
with the associative even part~%
$A$
and the odd part~%
$M$
which is an~associative~%
$A$-module
were described in \cite{Zh02,ShestZhel}.
The study in \cite{Zh02}
was considerably influenced by \cite{Sh98},
which described the simple $(-1,1)$-superalgebras of characteristic
  $\neq 2,3$.
In the Jordan case,
if a~superalgebra is not
the superalgebra of a~nondegenerate bilinear superform,
then its even part~%
$A$
is a~differentially simple algebra
with respect to some set of derivations,
and its odd part~%
$M$
is a~finitely generated projective~%
$A$-module of rank~1.
Here,
as for
$(-1,1)$-superalgebras,
we define multiplication in~%
$M$
using fixed finite sets of derivations and elements of~%
$A$.
It turns out that
every Jordan superalgebra of this type
is a~subsuperalgebra of the superalgebra of vector type
$J(\Gamma,D)$.
Under certain restrictions on~%
$A$
the odd part~%
$M$
is a~cyclic~%
$A$-module,
and consequently,
the original Jordan superalgebra
is isomorphic to the superalgebra
$J(\Gamma,D)$.
For instance,
if
$A$
is a~local algebra
then by the well-known Kaplansky theorem~%
$M$
is free,
and consequently,
it is a~cyclic~%
$A$-module.
If the ground field is of characteristic
$p>2$
then \cite{ShY} implies that~%
$A$
is a~local algebra;
thus,
$M$
is a~cyclic~%
$A$-module.
If
$A$
is the ring of polynomials in finitely many variables
then~%
$M$
is free  by \cite{Susl77},
and consequently,
it is a~cyclic~%
$A$-module.

A~natural question arose:
is the original superalgebra isomorphic to 
$J(\Gamma,D)$?
Equivalently,
is the odd part~%
$M$
a~cyclic~%
$A$-module?
Examples are constructed in
\cite{ShestZhel, Zhel09}
of unital simple special Jordan superalgebras
with certain associative even part
and the odd part~%
$M$
which is not free,
i.e.,
not cyclic.
In those examples the ground field
is either the field of real numbers
or an~arbitrary field of characteristic~0
in which the equation
$t^2+1=0$
has no solutions.

In this article we construct
a~similar example of a~Jordan superalgebra
 over an~arbitrary field of characteristic~0,
as well as an~example of a~simple Jordan superalgebra
which is a~subsuperalgebra of the Cheng--Kac Jordan superalgebra.
Examples of these superalgebras
answer a~question of Cantarini and Kac \cite{CantarKac}.

\medskip

Take a~field~%
$F$
of characteristic not equal to~2.
A~superalgebra
$J=J_0+J_1$
is a~%
${\rm Z_2}$-graded
$F$-algebra:
$$
J_0^2\subseteq J_0, J_1^2\subseteq J_0, J_1J_0\subseteq J_1,
J_0J_1\subseteq J_1.
$$
Put
$A=J_0$
and
$M=J_1$.
The spaces~%
$A$
and~%
$M$
are called the even and odd parts of~
$J$.
The elements of
$A\cup M$
are called homogeneous.
The expression
$p(x)$
with
$x\in A\cup M$
means the parity of~
$x$:
$p(x)=0$
for
$x\in A$
($x$
is even)
and
$p(x)=1$
for
$x\in M$
($x$
is odd).

Given~%
$x$
in~%
$J$
denote by
$R_x$
the operator of right multiplication by~%
$x$.
A~superalgebra~%
$J$
is called a~Jordan superalgebra
if the homogeneous elements satisfy
the operator identities
\begin{eqnarray}
& & aR_b=(-1)^{p(a)p(b)}bR_a,\label{scomm}
\end{eqnarray}
\begin{eqnarray}
& & R_{a^2}R_a=R_aR_{a^2}, \label{Jord}
\end{eqnarray}
\begin{eqnarray}
& &R_aR_bR_c+(-1)^{p(a)p(b)+p(a)p(c)+p(b)p(c)}R_cR_bR_a+
(-1)^{p(b)p(c)}R_{(ac)b}=\nonumber
\\ &
&\phantom{R_b}R_aR_{bc}+(-1)^{p(a)p(b)}R_bR_{ac}+(-1)^{p(a)p(c)+
p(b)p(c)}R_cR_{ab}. \label{supjor}
\end{eqnarray}

In every Jordan superalgebra,
the homogeneous elements satisfy
\begin{eqnarray}
&&(x,tz,y)=(-1)^{p(x)p(t)}t(x,z,y)+(-1)^{p(y)p(z)}(x,t,y)z,\label{supjor1}
\end{eqnarray}
where
$(x,z,y)=(xz)y-x(zy )$
is the associator of~
$x$,
$z$,
and~%
$y$.

Let us give some examples of Jordan superalgebras.

Take an~associative
$Z_2$-graded algebra
$B=B_0+B_1$
with multiplication~%
$\ast$.
Defining on the space~%
$B$
the supersymmetric product
$$
a\circ _s b=\frac{1}{2}(a\ast b+(-1)^{p(a)p(b)}b\ast a),
\mbox{ \ \ \ }
a,b\in B_0\cup B_1,
$$
we obtain the Jordan superalgebra
$B^{(+)s}$.
A~Jordan superalgebra
$J=A+M$
is called {\it special}
whenever it embeds
(as a~%
$Z_2$-graded algebra)
in the superalgebra
$B^{(+)s}$
for a~suitable
$Z_2$-graded associative algebra~%
$B$.

{\bf The superalgebra of vector type
$J(\Gamma,D)$.}
Take a~commutative associative~%
$F$-algebra~%
$\Gamma$
equipped with a~nonzero derivation~%
$D$.
Denote by
$\overline \Gamma$
an~isomorphic copy of the linear space~%
$\Gamma$,
and a fixed isomorphism,
by
$a\mapsto \overline a$.
On the direct sum
$J(\Gamma,D)=\Gamma+\overline \Gamma$
 of linear spaces
define a~multiplication~($\cdot$)
as
$$
a\cdot b=ab,
\mbox{ \ \ \ } a\cdot \overline b=\overline {ab},
\mbox{ \ \ } \overline a\cdot b=\overline {ab},
\mbox{ \ \ \ } \overline a\cdot
\overline b=D(a) b-aD(b),
$$
where
$a,b\in \Gamma$
and~%
$ab$
is the product in~%
$\Gamma$.
Then
$J(\Gamma,D)$
is a~Jordan superalgebra with the even part
$A=\Gamma$
and the odd part
$M=\overline \Gamma$.
The superalgebra
$J(\Gamma,D)$
is simple
if and only if
$\Gamma$
is a~%
$D$-simple algebra \cite{KMc90}
(i.e.,~%
$\Gamma$
contains no proper nonzero~%
$D$-invariant ideals,
and
$\Gamma^2=\Gamma$).

Consider the associative superalgebra
$B=M_2^{1, 1}({\rm End}\,\Gamma)$
with the even part
$$
B_0= \Big \{\left (
\begin{array}{cc}
  \phi & 0 \\
  0 & \psi
\end{array}\right),
\text { where } \phi,\psi\in {\rm End}\,\Gamma \Big \}
$$
and the odd part
$$
B_1= \Big\{\left (
\begin{array}{cc}
  0& \phi \\
  \psi & 0\end{array}\right)
\text { where }
\phi,\psi\in {\rm End}\,\Gamma\Big \}.
$$
It is shown in \cite{Mc} that
 the mapping
$$
a+\overline{b}\mapsto \left (
\begin{array}{cc}
 \text{\ \ } R_a & 4R_bD+2R_{D(b)} \\
  -R_b & R_a
\end{array}\right )
$$
is an~embedding of 
$J(\Gamma,D)$
into
$B^{(+)s}$.
Consequently,
the Jordan superalgebra
$J(\Gamma,D)$
is special.

{\bf The Kantor double
$J(\Gamma, \{\, ,\})$.}
Take an~associative supercommutative superalgebra
$\Gamma=\Gamma_0+\Gamma_1$
with unit~1
equipped with a~super-skew-symmetric bilinear mapping
$\{ \,, \}  : \Gamma\mapsto \Gamma$,
which we call the bracket.
From~
$\Gamma$
and~
$\{ \,, \}$
we can construct a~superalgebra
$J(\Gamma, \{\, ,\})$
as follows.
Consider the direct sum
$J(\Gamma, \{ \, ,\})=\Gamma\oplus \Gamma x$
of linear spaces,
where
$\Gamma x$
is an~isomorphic copy of~%
$\Gamma$.
Take two homogeneous elements~%
$a$
and~%
$b$
of~%
$\Gamma$.
The multiplication~%
($\cdot$)
on
$J(\Gamma, \{ \, ,\})$
is defined as
$$
a\cdot b=ab,
\mbox{ \ \ \ } a\cdot bx=(ab)x,
\mbox{ \ \ \ } ax\cdot b=(-1)^{p(b)}(ab)x,
\mbox{ \ \ \ } ax\cdot bx=(-1)^{p(b)}\{a,b\}.
$$
Put
$A=\Gamma_0+\Gamma_1x$
and
$M=\Gamma_1+\Gamma_0x$.
Then
$J(\Gamma, \{ \, ,\})=A+M$
is a
$Z_2$-graded algebra.

Refer to
$\{ \,,\}$
as a~Jordan bracket
if 
$J(\Gamma, \{ \, ,\})$
is a~Jordan superalgebra.
It is known (see \cite{KMc95}) that
  $\{ \,, \}$
is a~Jordan bracket
if and only if
it satisfies 
\begin{eqnarray}
& &
\{a,bc\}=\{a,b\}c+(-1)^{p(a)p(b)}b\{a,c\}-\{a,1\}bc,\label{skob1}
\end{eqnarray}
\begin{eqnarray}
& & \{a,\{b,c\}\}=\{\{a,b\},c\}+(-1)^{p(a)p(b)}\{b,\{a,c\}\}
+\{a,1\}\{b,c\}+ \nonumber \\ &
&(-1)^{p(a)(p(b)+p(c))}\{b,1\}\{c,a\}+(-1)^{p(c)(p(a)+p(b))}\{c,1\}\{a,b\},
\label{skob2}
\end{eqnarray}
\begin{eqnarray}
& & \{d,\{d,d\}\}=\{d,d\}\{d,1\}\label{skob3},
\end{eqnarray}
where
$a,b,c\in \Gamma_0\cup \Gamma_1$,
and
$d\in \Gamma_1$.

In particular,
$J(\Gamma,D)$
is the algebra
$J(\Gamma, \{\, ,\})$
if
$$\{a,b\}=D(a) b-aD(b).
$$

The next theorem is proved in \cite{ShestZhel}.
\medskip

{\bf Theorem.}
{\it
Take a~simple special unital Jordan superalgebra
$J=A+M$
whose even part~%
$A$
is an~associative algebra,
and whose odd part~%
$M$
is an~associative~%
$A$-module.
If~%
$J$
is not the superalgebra of a~nondegenerate bilinear superform
then there exist
$x_1,\ldots,x_n\in M$
such that
$$
M=x_1A+\ldots +x_nA,
$$
and the product in~%
$M$
satisfies
\begin{eqnarray}
& &ax_i\cdot bx_j=\gamma_{ij} ab+D_{ij}(a)b-aD_{ji}(b),\mbox{ \ \
} i,j=1,\ldots ,n,\label{Theor2}
\end{eqnarray}
where
$\gamma_{ij}\in A$,
and
$D_{ij}$
is a~derivation of~
$A$.
The algebra~%
$A$
is differentially simple with respect to the set of derivations
$\Delta \{D_{ij}|i,j=1,\ldots, n \}$.
The module~%
$M$
is a~projective~%
$A$-module of rank~1.
Moreover,
$J$
is a~subalgebra of the superalgebra
$J(\Gamma,D)$.  }
\medskip

In addition,
\cite{ShestZhel} includes
an~example of a~Jordan superalgebra
over the field of real numbers
satisfying the hypotheses of the theorem
which is not isomorphic to 
$J(\Gamma,D)$.
A~similar example of a~Jordan superalgebra
over a~field of characteristic zero
in which the equation
$t^2+1=0$
has no solutions
is constructed in \cite{Zhel09}.
Let us give another example of this kind of superalgebra
over an~arbitrary field of characteristic zero.

Fix an~arbitrary field~%
$F$
of characteristic~0.
Consider the polynomial algebra
$F[x,y]$
in two variables~%
$x$
and~%
$y$.
Denote by
$\partx$
and
$\party$
the operators of differentiation
with respect to~%
$x$
and~%
$y$
on
$F[x,y]$.
Put
$D=2y^3\partx-x\party$
and
$f(x,y)=x^2+y^4-1$.
Then~%
$D$
is a~derivation of
$F[x,y]$,
and
$D(f(x,y))=0$.
Take the quotient algebra
$\Gamma=F[x,y]/f(x,y)F[x,y]$
of
$F[x,y]$
by the ideal
$f(x,y)F[x,y]$.
It is clear that~%
$D$
induces a~derivation of~%
$\Gamma$,
which we denote by~%
$D$
as well.
Identify the images of~%
$x$
and~%
$y$
under the canonical homomorphism
$F[x,y]\mapsto \Gamma$
with the elements~%
$x$
and~%
$y$.
Then
$\Gamma=F[y]+xF[y]$,
where
$F[y]$
is the polynomial ring in~%
$y$.\medskip

{\bf Proposition 1.}
{\it The algebra~%
$\Gamma$
is differentially simple with respect to~%
$D$.}

{\sc Proof.}
Suppose that~%
$I$
is a~nonzero~%
$D$-invariant ideal of~%
$\Gamma$.
If
$f(y)\in F[y]$
and
$f(y)\in I$
then
$D(f(y))=-xf'(y)\in I$,
where
$f'(y)$
is the derivative of
$f(y)$
with respect to~%
$y$.
Then
$(1-y^4)f'(y)\in I$
and
$D((1-y^4)f'(y))\in I$.
Thus,
$$
-x(-4y^3f'(y)+(1-y^4)f''(y))\in I.
$$
This implies that
$(1-y^4)^2f''(y)\in I$.
Continuing this process,
we deduce that
$(1-y^4)^kf^{(k)}(y)\in I$
for all~%
$k$,
where
$f^{(k)}(y)$
is the order~%
$k$
derivative of
$f(y)$.
Consequently,
$(1-y^4)^k\in I$
for some~%
$k$.
Take the smallest~%
$k$
with
$z_k=(1-y^4)^k\in I$.
Then
$$
D(z_k)=4kxy^3(1-y^4)^{k-1}\in I.
$$
Thus,
$$
x(1-y^4)^{k-1}=xz_k+\frac{1}{4k}yD(z_k)\in I.
$$
Consequently,
$$
D(x(1-y^4)^{k-1})=2y^3(1-y^4)^{k-1}+(k-1)4y^3(1-y^4)^{k-1}2(2k-1)y^3(1-y^4)^{k-1}\in I.
$$
This implies that
$y^3(1-y^4)^{k-1}\in I$
and
$y^4(1-y^4)^{k-1}\in I$.
Then,
$$
z_{k-1}=(1-y^4)^k+y^4(1-y^4)^{k-1}\in I.
$$
Therefore,
we may assume that
$F[y]\cap I=0$.

Suppose that
$f(y)+xg(y)\in I$.
Then
$$
(f(y)+xg(y))(f(y)-xg(y))=f(y)^2-(1-y^4)g(y)^2\in I.
$$
By the argument above,
$f(y)^2=(1-y^4)g(y)^2$.
Then,
$1-y^4=h(y)^2$
for some 
$h(y)\in F[y]$,
and we arrive at a~contradiction.

Consequently,
$\Gamma$
is a~differentially simple algebra with respect to~%
$D$.\hfill
$\Box$
\medskip

Consider in~%
$\Gamma$
the subalgebra~%
$A$
generated by~%
$1$,
$y^2$,
and
$xy$.
Then,
$$
D(y^2)=-2xy\in A\text{ and  } D(xy)=3y^4-1\in A.
$$
Consequently,
$D(A)\subseteq A$.
Observe that
$1,y^{2i},xy^{2i-1}$,
where
$i=1,2,\dots$,
constitute a~linear basis for~%
$A$.
We can express every element of~%
$A$
as
$f(y)+xyg(y)$
with
$f(y),g(y)\in F[y^2]$.
 \medskip

{\bf Proposition 2.}
{\it The algebra~%
$A$
is differentially simple with respect to~%
$D$.}

{\sc Proof.}
Suppose that~%
$I$
is a~nonzero~%
$D$-invariant ideal of~%
$A$.
If
$f(y)\in F[y^2]$
and
$f(y)\in I$
then
$xf'(y)=-D(f(y))\in I$.
Thus,
$(1-y^4)yf'(y)=(xy)(xf'(y))\in I$.
Since
$$
D(xf'(y))=2y^3f'(y)-(1-y^4)f''(y)\in I,
$$
it follows that
$(1-y^4)^2f''(y)\in I$.
An~easy induction implies that
$$
(1-y^4)^{2k-1}yf^{(2k-1)}(y)\in I
\text{\ \  and \ \ }
(1-y^4)^{2k}f^{(2k)}(y)\in I.
$$
This yields
$(1-y^4)^{2k}\in I$.

Take the smallest~%
$k$
with
$(1-y^4)^{k}\in I$.
Then,
$$
D((1-y^4)^{k})=-4kxy^3(1-y^4)^{k-1}\in I.
$$
Consequently,
$$
xy(1-y^4)^{k-1}=xy(1-y^4)^{k}+y^2(xy^3(1-y^4)^{k-1})\in I.
$$
Thus,
$$
D(xy(1-y^4)^{k-1})=(3y^4-1)(1-y^4)^{k-1}+(k-1)4y^4(1-y^4)^{k-1}((4k-1)y^4-1)(1-y^4)^{k-1}\in I.
$$
Then,
$$
(4k-2)(1-y^4)^{k-1}=(4k-1)(1-y^4)^{k}+((4k-1)y^4-1)(1-y^4)^{k-1}\in I.
$$
Therefore,
we may assume that
$F[y^2]\cap I=0$.

Suppose that
$f(y)+xyg(y)\in I$.
Then,
$$
f(y)^2-(1-y^4)y^2g(y)^2=(f(y)+xyg(y))(f(y)-xyg(y))\in I.
$$
By the argument above,
$f(y)^2-(1-y^4)y^2g(y)^2=0$,
and we arrive at a~contradiction
since
$\deg f(y)^2=4n$
but
$\deg (1-y^4)y^2g(y)^2=4m+6$.

Therefore,
$A$
is a~differentially simple algebra with respect to~%
$D$.\hfill
$\Box$\medskip

The subspace
$M=xA+yA$
of~%
$\Gamma$
is an~associative~%
$A$-module.
\medskip

{\bf Proposition 3.}
{\it The module~%
$M$
is not a~cyclic~%
$A$-module.}

{\sc Proof.}
Assuming the contrary,
denote the generator of~%
$M$
 by~%
$z$.
Then
$z=xa+yb$
with
$a,b\in A$,
$x=zc$,
and
$y=zd$
for some
$c,d\in A$.
This implies that
\begin{equation}\label{Eq1}
xd=yc,
\end{equation}
\begin{equation}\label{Eq2}
x=x(ac+bd), y=y(ac+bd).
\end{equation}
We can write
$$
a=f_0+xyf_1, b=g_0+xyg_1, c=e_0+xye_1, d=h_0+xyh_1,
$$
where
$f_0,f_1,g_0,g_1,e_0,e_1,h_0,h_1$
are polynomials in
$F[y^2]$.

From (\ref{Eq1}) we deduce that
$$
h_0=y^2e_1\text{ \ \ and \ \ }
e_0=(1-y^4)h_1.
$$
From  (\ref{Eq2}) we deduce that
\begin{equation}\label{Eq3}
f_0e_0+(1-y^4)y^2f_1e_1+g_0h_0+(1-y^4)y^2g_1h_1=1,
\end{equation}
\begin{equation}\label{Eq4}
f_0e_1+f_1e_0+g_0h_1+g_1h_0=0.
\end{equation}

Denote by
$(e_1,h_1)$
the greatest common divisor of 
$e_1$
and
$h_1$.
Since
$h_0=y^2e_1$
and
$e_0=(1-y^4)h_1$,
by (\ref{Eq3}) we have
$$
1=(1-y^4)f_0h_1+(1-y^4)y^2f_1e_1+y^2g_0e_1+(1-y^4)y^2g_1h_1=$$
$$(1-y^4)(f_0+y^2g_1)h_1+y^2((1-y^4)f_1+g_0)e_1.
$$
Consequently,
$(e_1,h_1)=1$.
By (\ref{Eq4}),
\begin{equation}
(f_0+y^2g_1)e_1+((1-y^4)f_1+g_0)h_1=0.\nonumber
\end{equation}
This and
$(e_1,h_1)=1$
imply that
$f_0+y^2g_1=h_1u$,
where
$u\in F[y]$.
Then,
$$
h_1ue_1+((1-y^4)f_1+g_0)h_1=0.
$$
Thus,
$$
ue_1+((1-y^4)f_1+g_0)=0.
$$
By the argument above,
$$
1=(1-y^4)(f_0+y^2g_1)h_1+y^2((1-y^4)f_1+g_0)e_1=(1-y^4)h_1^2u-y^2e_1^2u.
$$
Then,
$u\in F$.
Consequently,
$$
(1-y^4)h_1^2u=1+y^2e_1^2u,
$$
which is impossible
since on the left we have a~polynomial of degree
$4k+4$,
while on the right,
of degree
$4m+2$.

Therefore,
$M$
is not a~cyclic~%
$A$-module.\hfill
$\Box$\medskip

Put
$$
D_{11}=(1-y^4)D,D_{12}= xyD, D_{22}=y^2D.
$$
Then
$D_{11},D_{12}, D_{22}$
are derivations of~
$A$.\medskip

{\bf Proposition 4.}
{\it The algebra~%
$A$
is differentially simple with respect to the set of derivations
$\Delta=\{D_{11}, D_{12}, D_{22}\}$.}

{\sc Proof.}
Suppose that
$I$
is an~ideal of~
$A$
closed under 
$\Delta$. 
Then
$y^2D_{22}(I)\subseteq y^2I\subseteq I$.
Since
$$
D=D_{11}+y^2D_{22},
$$
it follows that
$D(I)\subseteq I$.
By Proposition~2,
either
$I=0$
or
$I=A$.
Consequently,~%
$A$
is a~differentially simple algebra with respect to
$\Delta=\{D_{11},D_{12}, D_{22}\}$.\hfill
$\Box$\medskip

Consider now the superalgebra
$J(\Gamma,D)$.
Proposition~1 implies that
$J(\Gamma,D)$
is a~simple superalgebra.
Consider its subspace
$$
J(A,\Delta)=A+\overline{M}.
$$
Recall that~%
$A$
is the subalgebra of~%
$\Gamma$
generated by~
$1$,
$y^2$,
and
$xy$,
while
$M=xA+yA$.

Given
$a,b\in A$,
in 
$J(\Gamma,D)$
we have
$$
\overline{xa}\cdot \overline{xb}=D(xa)xb-D(xb)xa=$$
$$D(x)axb+D(a)x^2b-D(x)xab-D(b)x^2a=
D_{11}(a)b-aD_{11}(b)\in A.
$$
Similarly,
$$
\overline{ya}\cdot
\overline{yb}=D(y)ayb+D(a)y^2b-D(y)yab-D(b)y^2a=D_{22}(a)b-aD_{22}(b)\in A,
$$
$$
\overline{xa}\cdot
\overline{yb}=D(x)ayb+D(a)xyb-D(y)xab-D(b)yxa=(1+y^4)ab+D_{12}(a)b-aD_{12}(b)\in A.
$$
Consequently,
$J(A,\Delta)$
is a~subsuperalgebra of
$J(\Gamma,D)$.
Thus,
$J(A,\Delta)$
is a~Jordan superalgebra.
Moreover,
the odd elements in
$J(\Gamma,D)$
multiply according to (\ref{Theor2}),
where
$\Delta=\{D_{11},D_{12}, D_{22}\}$,
and
$\gamma_{12}=1+y^4$.
By Proposition~3,
$J(A,\Delta)$
is not isomorphic to a~superalgebra of type
$J(\Gamma_0,D_0)$.

Verify that
$J(A,\Delta)$
is a~simple superalgebra.
Suppose that~%
$I$
is a~nonzero
${\rm Z_2}$-graded ideal of
$J(A,\Delta)$.
Then
$I=I_0+I_1$,
where
$I_0$
is an~ideal of~%
$A$.
Given
$r\in I_0$,
we have
$$D_{11}(r)=\overline{(xr)}\cdot
\overline{x}=(r\cdot\overline{x})\cdot \overline{x} \in I_0.$$
Similarly, $D_{12}(r), D_{22}(r)\in I_0$. Consequently, $I_0$
is invariant under the set of derivations~%
$\Delta$.
By Proposition~4,
  either
$I_0=A$
or
$I_0=0$.
If
$I_0=A$
then
$1\in I_0\subseteq I$
and
$I=J(A,\Delta)$.
If
$I_0=0$
then
$I\subseteq \overline{M}$
and
$I\cdot \overline{M}\subseteq I_0=0$.
It is clear that
$$
A=AD_{11}(A)+AD_{12}(A)+AD_{22}(A).
$$
Thus,
$$
1=\sum_i(a_{1i},\overline{x},\overline{x})b_{1i}+\sum_i(a_{2i},\overline{x},\overline{y})b_{2i}+
\sum_i(a_{3i},\overline{y},\overline{y})b_{3i}
$$
for some elements
$a_{1i}$,
$a_{2i}$,
$a_{3i}$,
$b_{1i}$,
$b_{2i}$,
and
$b_{3i}$
of
$A$.
By (\ref{supjor1}) we deduce that
$1\in (A,\overline{M},\overline{M})$
and
$$
I\cdot(A,\overline{M},\overline{M})\subseteq (A,I\cdot
\overline{M},\overline{M})+(A,I,\overline{M})\cdot
\overline{M}=0.
$$
Then,
$I=0$.
Consequently,
$J(A,\Delta)$
is a~simple superalgebra.

Let us summarize the argument as
\medskip

{\bf Theorem 1.}
{\it
Take an~arbitrary field~%
$F$
of characteristic~0.
Consider the polynomial algebra
$F[x,y]$
in two variables~%
$x$
and~%
$y$.
Put
$f(x,y)=x^2+y^4-1$
and
$D=2y^3\partx-x\party$.
Put
$\Gamma=F[x,y]/f(x,y)F[x,y]$.
Then the derivation~%
$D$
induces a~derivation of the algebra~%
$\Gamma$,
which we denote by~%
$D$
as well.
Identify the images of~%
$x$
and~%
$y$
under the canonical homomorphism
$F[x,y]\mapsto \Gamma$
with the elements~%
$x$
and~%
$y$.
Suppose that~%
$A$
is a~subalgebra of~%
$\Gamma$
generated by~%
$1$,
$y^2$,
and
$xy$,
while
$M=xA+yA$.
Put
$$
\Delta=\{D_{11},D_{12}, D_{22}\},
\text{ where }
D_{11}=(1-y^4)D,D_{12}= xyD, D_{22}=y^2D.
$$
Then the subspace
$J(A,\Delta)=A+\overline{M}$
is a~subsuperalgebra of
$J(\Gamma,D)$,
and the multiplication of odd elements in
$J(A,\Delta)$
is defined as
$$
\overline{xa}\cdot\overline{xb}=D_{11}(a)b-aD_{11}(b),
\text{ \ \ }
\overline{ya}\cdot\overline{yb}=D_{22}(a)b-aD_{22}(b),
$$
$$
\overline{xa}\cdot\overline{yb}=(1+y^4)ab+D_{12}(a)b-aD_{12}(b).
$$
Moreover,
$J(A,\Delta)$
is a~simple superalgebra,
and~%
$\overline{M}$
is not a~cyclic~%
$A$-module;
i.e.,
$J(A,\Delta)$
is not isomorphic to a~superalgebra of vector type
$J(\Gamma_0,D_0)$.}
\medskip

{\bf The Superalgebra of Type
$JS(\Gamma,D)$.}
Take an~associative supercommutative superalgebra
$\Gamma=\Gamma_0+\Gamma_1$
equipped with a~nonzero odd derivation~%
$D$;
i.e.,
$D(\Gamma_i)\subseteq \Gamma_{(i+1)mod\, 2}$
and
$$D(ab)=D(a)b+(-1)^{p(a)}aD(b)$$
for
$a,b\in\Gamma_0\cup\Gamma_1$.

Put
$A=\Gamma_1$,
$M=\Gamma_0$,
and
$JS(\Gamma,D)=A+M$.
Define on the space
$JS(\Gamma,D)$
the multiplication
$$
a\circ b=aD(b)+(-1)^{p(a)}D(a)b.
$$
Then
$JS(\Gamma,D)$
is a~Jordan superalgebra.
If
$JS(\Gamma,D)$
is a~simple superalgebra
then~%
$\Gamma$
is a~differentially simple superalgebra
(see~\cite{CantarKac}).
\medskip

{\bf Proposition 5.}
{\it The superalgebra
$JS(\Gamma,D)$
is not unital.}

{\sc Proof.}
Suppose that
$e$
is the unit of 
$JS(\Gamma,D)$.
Then
$e\in A\subseteq \Gamma_1$.
Given
$a\in JS(\Gamma,D)$,
we have
$$
a=e\circ a=eD(a)+D(e)a.
$$
Since~%
$\Gamma$
is supercommutative
and
$e\in \Gamma_1$,
it follows that
$e=2eD(e)$
and
$e^2=0$
in~%
$\Gamma$.
Consequently,
$ea=eD(e)a=\frac{1}{2}ea$.
This implies that
$e\Gamma=0$.
Then,
$e=2eD(e)=0$.\hfill
$\Box$\medskip

{\bf Corollary 1.}
{\it The superalgebra
$J(A,\Delta)$
is not isomorphic to the superalgebra
$JS(\Gamma,D)$.}\medskip

{\bf The Cheng--Kac superalgebra.}
Take an~associative commutative~%
$F$-algebra~%
$\Gamma$
equipped with a~nonzero derivation~%
$D$.
Consider two direct sums
$$
J_0=\Gamma+w_1\Gamma+w_2\Gamma+w_3\Gamma
$$
and
$$
J_1=\overline{\Gamma}+x_1\overline{\Gamma}+x_2\overline{\Gamma}+x_3\overline{\Gamma}
$$
of linear spaces,
where~%
$\overline{\Gamma}$
is an~isomorphic copy of~%
$\Gamma$.

For
$a,b\in \Gamma$
define a~multiplication on the space
$J_0$
by putting
$$
a\cdot b=ab,\, a\cdot w_ib=w_iab,\, w_1a\cdot w_1b=w_2a\cdot
w_2b=ab,\, w_3a\cdot w_3b=-ab,$$ $$w_ia\cdot w_jb=0 \text{ for }
i\neq j.
$$
Put
$x_{i\times i}=0$,
$x_{1\times 2}=-x_{2\times 1}=x_3$,
$x_{1\times 3}=-x_{3\times 1}=x_2$,
and
$x_{2\times 3}=-x_{3\times 2}=-x_1$.
Define a~bimodule action
$J_0\times J_1\mapsto J_1$
by putting
$$
a\cdot\overline{b}=\overline{ab},\,
a\cdot x_i\overline{b}=x_i\overline{ab},\,
w_ia\cdot\overline{b}=x_i\overline{D(a)b},\,
w_ia\cdot x_j\overline{b}=x_{i\times j}\overline{ab}.
$$
The bracket on
$J_1$
is defined as
$$
\overline{a}\cdot \overline{b}=D(a)b-aD(b),\,
\overline{a}\cdot x_i\overline{b}=-w_i(ab),\,
x_i\overline{a}\cdot \overline{b}=w_i(ab),\,
x_i\overline{a}\cdot x_j\overline{b}=0.
$$
Then the space
$J=J_0+J_1$
with the multiplication
$$
(a_0+a_1)\cdot(b_0+b_1)=(a_0\cdot b_0+a_1\cdot b_1)+(a_0\cdot
b_1+a_1\cdot b_0)
$$
for
$a_0, b_0\in J_0$
and
$a_1, b_1\in J_1$
is an~algebra,
which is denoted by
$CK(\Gamma, D)$.
It is known
(see \cite{MZ2001, CantarKac})
that
$CK(\Gamma, D)$
is a~Jordan superalgebra,
which is simple
if and only if~%
$\Gamma$
is~%
$D$-simple.

Suppose now that $\Gamma=F[x,y]/f(x,y)F[x,y]$, where
$f(x,y)=x^2+y^4-1$ and $D=2y^3\partx-x\party$. Consider the Jordan
superalgebra $J(A,\Delta)=A+\overline{M}$ constructed above.
In 
$CK(\Gamma,D)$
consider the subspace
$$
GCK(A,\Delta)=A+w_1A+w_2A+w_3A+\overline{M}+x_1\overline{M}+x_2\overline{M}+x_3\overline{M}.
$$
In~
$\Gamma$
we have
$M^2\subseteq A$.
Thus,
$GCK(A,\Delta)$
is a~subsuperalgebra of
$CK(\Gamma,D)$.
Consequently,
$GCK(A,\Delta)$
is a~Jordan superalgebra with the even part
$GCK(A,\Delta)_0=A+w_1A+w_2A+w_3A$
and the odd part
$GCK(A,\Delta)_1=\overline{M}+x_1\overline{M}+x_2\overline{M}+x_3\overline{M}$.
\medskip

{\bf Theorem 2.}
{\it
For an~arbitrary field~%
$F$
of characteristic zero
$GCK(A,\Delta)$
is a~simple unital Jordan superalgebra. }

{\sc Proof.}
Suppose that
$I=I_0+I_1$
is a~nonzero ideal of 
$GCK(A,\Delta)$.
Then
$K=A\cap I_0$
is an~ideal of~
$A$,
and
$(K,\overline{M},\overline{M})\subseteq K$.
Thus,
$K+K\cdot \overline{M}$
is an~ideal of 
$J(A,\Delta)$.
If
$K\neq 0$
then since
$J(A,\Delta)$
is a~simple superalgebra,
we have
$1\in K$.
Consequently,
$I=GCK(A,\Delta)$.

Suppose that
$A\cap I_0=0$
and take
$r=a+w_1a_1+w_2a_2+w_3a_3\in I_0$.
Then
$w_2(w_2(w_1r))=a_1\in A\cap I_0$.
Consequently,
$a_1=0$.
Similarly,
$a_2=a_3=0$.
Thus,
$I_0=0$.
This implies that
$I\subseteq GCK(A,\Delta)_1$
and
$I\cdot GCK(A,\Delta)_1\subseteq I_0=0$.
Since
$1\in (A,\overline{M},\overline{M})$,
by (\ref{supjor1}) we deduce that
$$
I\cdot(A,\overline{M},\overline{M})\subseteq (A,I\cdot
\overline{M},\overline{M})+(A,I,\overline{M})\cdot
\overline{M}=0.
$$
Then,
$I=0$.
Consequently,
$GCK(A,\Delta)$
is a~simple superalgebra.\hfill
$\Box$\bigskip

I would like to take this chance
to express by special gratitude to A.~P.~Pozhidaev,
whose comments helped to improve this article.

\noindent
ZHELYABIN Viktor Nikolaevich, \\
Sobolev Institute of Mathematics, RAS\\
4 Acad. Koptyug prospekt\\
Novosibirsk 630090\\
RUSSIA\\
\smallskip
phone +7(383)(363-45-57)\\
email: {\tt vicnic@math.nsc.ru}\\
\smallskip
and\\
Novosibirsk State University\\
2 Pirogova str.\\

\end{document}